\documentclass[11pt]{article}
\usepackage{amsmath}
\usepackage{amssymb}
\usepackage{graphicx}

\newtheorem{theorem}{Theorem}%[section]

\def\cE{{\mathcal E}}
\def\qed{\hfill $\Box$\medskip}

\def\IC{{\mathbb C}}

\def\IR{{\mathbb R}}

\def\tr{{\mathrm tr}\,}

\def\diag{{\rm diag}\,}

\def\({\left(}
\def\){\right)}
\def\[{\left[}
\def\]{\right]}

\def\dfrac{\displaystyle\frac}

\begin{document}
\openup .38\jot
\title{Submultiplicativity of the numerical radius of
commuting matrices of order two}
\author{Chi-Kwong Li\footnote{
Department of Mathematics, College of William and Mary, Williamsburg, VA 23187. 
(ckli@math.wm.edu)}{\ }   and Yiu-Tung Poon\footnote{Department of Mathematics, Iowa Sate University,
Ames, IA 50011. (ytpoon@iastate.edu)}}
\date{}
\maketitle

\centerline{\bf Dedicated to Professor Pei Yuan Wu.}

\begin{abstract}
Denote by $w(T)$ the numerical radius of a matrix $T$.
An elementary proof is given to the fact
that $w(AB) \le w(A)w(B)$ for a pair of commuting
matrices of order two, and  characterization is given 
for the matrix pairs that attain the quality.
\end{abstract}



\section{Introduction}

Let $M_n$ be the set of $n\times n$ matrices. The numerical range and numerical radius
of $A \in M_n$ are defined by
$$W(A) = \{x^*Ax: x \in \IC^n, x^*x = 1\} \qquad \hbox{ and } \qquad
w(A) = \max\{|\mu|: \mu \in W(A)\},$$
respectively.
The numerical range and numerical radius are useful tools in studying matrices and
operators.
There are strong connection between the algebraic properties
of a matrix $A$ and the geometric properties of $W(A)$. 
For example, $W(A) = \{\mu I\}$ if and only if $A = \mu I$;
$W(A) \subseteq \IR$ if and only if $A = A^*$; $W(A) \subseteq [0, \infty)$
if and only if $A$ is positive semi-definite.

The numerical radius is a norm on $M_n$, and 
has been used in the analysis of basic iterative solution  
methods  \cite{Ax}.  Researchers have 
obtained interesting inequalities related to the numerical radius;
for example, see
\cite{G,H,Hol,HolS,HJ1}.  We mention a few of them in the following.
Let $\|A\|$
% = \max\{ (x^*A^*Ax)^{1/2}: x \in \IC^n, x^*x = 1 \}$
be the operator norm of $A$.
It is known that
$$w(A) \le \|A\| \le 2w(A).$$
While the spectral norm is submultiplicative,
i.e., $\|AB\| \le \|A\| \|B\|$ for all $A, B \in M_n$,
the numerical radius is not.
In  general,
$$w(AB) \le \xi w(A)w(B) \quad \hbox{ for all } A, B \in M_n$$
if and only if $\xi \ge 4$; e.g., see  \cite{GW}.
Despite the fact that the numerical radius is not submultiplicative,
$$w(A^m) \le w(A)^m \qquad \hbox{ for all positive integers } m.$$
For a normal matrix $A\in M_n$, we have $w(A) = \|A\|$. Thus, 
for a normal matrix $A$ and any $B \in M_n$,
$$w(AB) \le \|AB\| \le \|A\| \|B\| = w(A) \|B\| \le 2 w(A)w(B),$$
and also
$$w(BA) \le \|BA\| \le \|B\|\|A\|  = \|B\| w(A) \le 2w(B)w(A).$$
In case  $A, B\in M_n$ are normal matrices,
$$w(AB) \le \|AB\| \le \|A\| \|B\| = w(A) w(B).$$
Also, for any pairs of commuting matrices $A, B \in M_n$,
$$w(AB) \le 2w(A)w(B).$$
To see this, we may assume $w(A) = w(B) = 1$, and observe that
\begin{eqnarray*}
4w(AB) &=&  w((A+B)^2 - (A-B)^2) \le w((A+B)^2) + w((A-B)^2) \\
&\le& w(A+B)^2 + w(A-B)^2 \le 8.
\end{eqnarray*}  
The constant 2 is best (smallest) possible for matrices of order at least 4 
because $w(AB) = 2w(A)w(B)$ if $A = E_{12} + E_{34}$ and $B = E_{13} + E_{24}$, where $ E_{ij} \in M_n$ has $1$ at the $(i,j)$ position and $0$ elsewhere; see \cite[Theorem 3.1]{GW}.

In connection to the above discussion, 
there has been interest in studying the best (smallest) constant
$\xi > 0$ such that
$$w(AB) \le \xi w(A)w(B)$$
for all commuting matrices $A, B \in M_n$ with $n \le 3$.
For $n = 2$, the best constant $\xi$ is one;
the existing proofs of the $2\times 2$ case depend on deep theory on analytic functions,
von Neumann inequality, and functional calculus
on operators with numerical radius equal to one, etc.; for example, see \cite{Hol,HolS}.

Researchers have been trying to find an elementary proof for this result in view of the fact
that the numerical range of $A \in M_2$ is well understood, 
namely,  $W(A)$is an elliptical disk with the eigenvalues
$\lambda_1, \lambda_2$ as
foci and the length of minor axis $\sqrt{(\tr A^*A) - |\lambda_1|^2 - |\lambda_2|^2}$;
for example, see \cite{L,M} and \cite[Theorem 1.3.6]{HJ1}.

The purpose of this note is to provide such a  proof.
Our analysis is based on elementary theory in convex analysis, co-ordinate
geometry, and inequalities. 
Using our approach, we readily  give a  characterization of commuting pairs of matrices  
$A, B \in M_2$ satisfying $w(AB) = w(A)w(B)$, which was done
in \cite[Theorem 4.1]{GW} using yet another deep result of Ando \cite{An} that  a
matrix $ A$ has numerical radius bounded by one if and only if $A = (I-Z)^{1/2}C(A+Z)^{1/2}$
for some contractions $C$ and $Z$, where $Z = Z^*$.
Here is our main result.

\begin{theorem} Let $A, B \in M_2$ be nonzero matrices such that $AB = BA$.
Then $w(AB) \le w(A)w(B)$. The equality holds if and only if one of the following holds.

\begin{itemize}
\item[{\rm (a)}] $A$ or $B$ is a scalar matrix, i.e. of the form $ \mu I_2$ for some 
$\mu \in\IC$.
\item[
{\rm (b)}] There is a unitary $U$ such that $U^*AU = \diag(a_1,a_2)$
and $U^*BU = \diag(b_1, b_2)$ with $|a_1| \ge |a_2|$ and $|b_1| \ge |b_2|$.
\iffalse
\item[
{\rm (c)}] There is a unitary $U$ such that
$U^*AU = a R$ and $ U^*BU =b\overline R$ for some $a, b \in \IC$ with $|a|= w(A)$, $|b| = w(B)$,
$R =  ic^2 \gamma I_2 + \sqrt{1- (c\gamma )^2}
\begin{pmatrix} c & 2\sqrt{1-c^2} \cr 0
& -c \end{pmatrix}$ for some $\gamma, c \in (0,1]$.
\fi
\end{itemize}
\end{theorem}

One can associate the conditions (a) and (b) in the theorem with the geometry of the numerical
range of $A$ and $B$ as follows.
Condition (a) means that $W(A)$ or $W(B)$ is a single point; condition (b) means that
$W(A)$, $W(B)$, $W(AB)$ are line segments with three sets of
end points,
$\{a_1, a_2\}, \{b_1, b_2\}, \{a_1b_1, a_2b_2\}$, respectively, such that
$|a_1|\ge |a_2|$ and $|b_1| \ge |b_2|$. 

\iffalse
condition (c) means that there are $a, b \in \IC$ with $(|a|, |b|) = w(A), w(B))$ such that
$W(A/a) = \cE$ is an elliptical disk symmetric about the imaginary axis insider the unit 
circle at two points, and $W(B/b) = \overline{\cE}   =
\{\bar\mu \in \IC: \mu \in \cE\}.$
\fi

\section{Proof of Theorem 1}
Let $A, B \in M_2$ be commuting matrices. We may replace
$(A,B)$ by $(A/w(A),B/w(B))$ and assume that $w(A) = w(B) = 1$.
We need to show that $w(AB)\le 1$.

Since $AB = BA$, there is a unitary matrix $U \in M_2$ such that
both $U^*AU$ and$U^*BU$ are in triangular form; for example, see 
\cite[Theorem 2.3.3]{HJ2}. 
We may replace $(A,B)$ by $(U^*AU,U^*BU)$ and assume that
$A = \begin{pmatrix} a_1&a_3\cr 0 & a_2\cr \end{pmatrix}$,
$B = \begin{pmatrix} b_1&b_3\cr 0 & b_2\cr \end{pmatrix}$  and  $w(A) = w(B) = 1$.
The result is clear if $A$ or $B$ is normal. So, we assume that $a_3, b_3 \ne 0$.
Furthermore, comparing the $(1,2)$ entries on both sides of $AB = BA$, we see that
$\dfrac{a_1-a_2}{a_3}= \dfrac{b_1-b_2}{b_3}$.  
Applying a diagonal unitary similarity to both $A$ and $B$, we may further assume that $\gamma =\dfrac{a_1-a_2}{a_3}\ge 0$. Let   
$r=\dfrac{1}{\sqrt{\gamma^2+1}}$. We have $0<r\le 1$. Then
$A = z_1 I + s_1C$ and $B = z_2 I + s_2 C$ with 

\bigskip
\centerline{
$z_1 = \dfrac{a_1+a_2}{2}, \quad
z_2 =\dfrac{b_1+b_2}{2} , \quad s_1=\dfrac{a_3}{2r}, \quad s_2  =\dfrac{b_3}{2r}$, 
\quad and \quad
$C = \begin{pmatrix} \sqrt{1-r^2} & 2r \cr 0 & -\sqrt{1-r^2}\end{pmatrix}$.}

\medskip\noindent
Note that $W(C)$ is the elliptical disk with boundary
$$\{\cos\theta + i r\sin\theta: \theta \in [0, 2\pi]\};$$ 
see \cite {L} and \cite[Theorem 1.3.6]{HJ1}.
Replacing $(A, B)$ with $(e^{it_1}  A, e^{it_2 }  B)$
for suitable $t_1, t_2\in  [0,2\pi ]$, if necessary, 
we may assume that
${\rm Re}\,z_1,\ {\rm Re}\,z_2 \ge 0 $ and $s_1,s_2$ are real.

Suppose   $z_1 = \alpha_1 + i \alpha_2$ with $\alpha_1\ge 0$  and the boundary of $W(A)$ touches the unit 
circle at the point $\cos\phi_1+i\sin\phi_1$ with $\phi_1 \in
[-\pi/2, \pi/2]$. 
Then $W(A)$ has boundary 
$$\{\alpha_1 +| s_1|\cos\theta + i (\alpha_2 +  |s_1|r\sin\theta): \theta \in [0, 2\pi]\} .$$ 

\noindent

We {\bf claim} that the matrix $A$ is a convex combination of
$A_0 = e^{i\phi_1} I$ and another matrix $A_1$ of the form
$A_1 = i(1-r^2)\sin \phi_1I + \xi C$ for some $\xi \in \IR$ such that $w(A_1) \le 1$.

To prove our claim, 
we first determine $\theta_1 \in [-\pi/2, \pi/2]$   satisfying
$$\cos\phi_1+i\sin\phi_1=(\alpha_1 + |s_1| \cos\theta_1) + i (\alpha_2 
+ |s_1| r\sin\theta_1).$$
Since the boundary of $W(A)$ touches the unit 
circle at the point $\cos\phi_1+i\sin\phi_1$, using the parametric equation
\begin{equation}\label{para}
x+iy = (\alpha_1 + |s_1| \cos\theta) + i (\alpha_2  + |s_1| r\sin\theta),
\end{equation}
 of the boundary of $W(A)$, 
we see that the direction of the tangent 
at the intersection point $\cos\phi_1+i\sin\phi_1$
is $-  \sin \theta_1 + i  r \cos \theta_1$, which agrees with $-\sin\phi_1+i \cos \phi_1 $, the  
direction of the  tangent line of the unit circle at the same point.
As a result, we have
$$(\cos \theta_1, \sin \theta_1)
= \dfrac{(\cos \phi_1, r\sin\phi_1)}{\sqrt{\cos^2\phi_1 + r^2\sin^2\phi_1}}.$$
Furthermore,   since
$\cos\phi_1+i\sin\phi_1=(\alpha_1 + |s_1| \cos\theta_1) 
+ i (\alpha_2  + |s_1| r\sin\theta_1)$,
we have
$$
\alpha_1=\cos\phi_1 - \dfrac{|s_1|\cos\phi_1}{\sqrt{\cos^2\phi_1 + r^2\sin^2\phi_1}}\ge 0 \quad
 \hbox{ and } \quad
\alpha_2=\sin \phi_1-\dfrac{|s_1|r^2\sin\phi_1}{\sqrt{\cos^2\phi_1 + r^2\sin^2\phi_1}}.
 $$
\noindent
\bf Assertion. \it If
 $\hat s_1 = \sqrt{\cos^2\phi_1+r^2 \sin^2 \phi_1}$, then  
$|s_1|\le \hat s_1$.
\rm

If $\cos \phi_1>0$, then 
$\alpha_1=\(1 - \dfrac{|s_1| }{\hat s_1}\)\cos\phi_1\ge 0$, and hence
$|s_1|\le \hat s_1$.

If $\cos\phi_1=0$, then $\sin\phi_1=\pm 1$, $\hat s_1=r$ and 
 $(\alpha_1,\alpha_2) = (0,\sin\phi_1 (1-|s_1|r))$ 
so that the parametric equation of
  the boundary of $W(A)$ in (\ref{para}) becomes 
$$x+iy =  |s_1| \cos\theta  
+ i ( \sin\phi_1 (1-|s_1|r)\ + |s_1| r\sin\theta)\,.$$
Since $w(A)=1$ and $\sin\phi_1=\pm 1$, for all $\theta\in [0,2\pi)$ ,  we have  
\begin{eqnarray*}
0 &\le&  1 - \[(|s_1| \cos\theta)^2+(\sin\phi_1 (1-|s_1|r)  
+ |s_1| r\sin\theta)^2 \]\\
&= &1- \[ |s_1| (1-\sin^2\theta) +(\pm (1-|s_1|r)   
+ |s_1| r\sin\theta)^2 \]\\
&= &1- \[ |s_1|^2  (1-(\pm 1\mp(1\mp\sin\theta))^2) +(   1-|s_1|r(1\mp    
   \sin\theta))^2 \]\\
&= &1- \[ |s_1|^2  (  2(1\mp\sin\theta)      - (1\mp\sin\theta)^2) +    1-2|s_1|r(1\mp    
   \sin\theta)+  |s_1|^2r^2(1\mp    
   \sin\theta)^2     \]\\
&= & 2|s_1|(r-|s_1|)(1\mp\sin\theta)+(1-r^2)|s_1|^2(1\mp    
   \sin\theta)^2  .
\end{eqnarray*}
Therefore, $(r-|s_1|)\ge 0$, which gives $|s_1|\le r=\hat s_1$.

\medskip
Now, we show that our {\bf claim} holds with  
\begin{equation}\label{a0a1}
A_0=e^{i\phi_1}I \qquad  \hbox{ and } \qquad  A_1 = i(1-r^2)\sin\phi_1I+\nu_1 \hat s_1 C,
\end{equation}
where  $\nu_1 =1$   if $s_1\ge 0$ and $\nu_1 =-1$ if $s_1< 0$.

Note that
$W(A_1)$ is the elliptical disk with boundary
$\{ \hat s_1 \cos \theta + i[(1-r^2) \sin\phi_1+ \hat s_1r\sin\theta): \theta \in [0, 2\pi)\}$,
and for every $\theta\in [0, 2\pi]$, we have 
\begin{eqnarray*}
&&(\hat s_1\cos \theta)^2+((1-r^2)\sin\phi_1+\hat s_1r\sin \theta)^2\\  
&=&\hat s_1^2(1-\sin^2 \theta)+(1-r^2)^2\sin^2\phi_1+\hat s_1^2r^2\sin^2 \theta+2\hat s_1r(1-r^2)\sin\phi_1\sin \theta\\ 
&=&\hat s_1^2   +(1-r^2)^2\sin^2\phi_1+(1-r^2)r^2\sin^2\phi_1    -(1-r^2)\(\hat s_1^2\sin^2 \theta-2\hat s_1r \sin\phi_1\sin \theta+r^2\sin^2\phi_1\)\\ 
&=&(\cos^2\phi_1+r^2 \sin^2 \phi_1) +(1-r^2)^2\sin^2\phi_1+(1-r^2)r^2\sin^2\phi_1-(1-r^2)(\hat s_1\sin \theta-r\sin\phi_1)^2
\\ 
&=&1-(1-r^2)(\hat s_1\sin \theta-r\sin\phi_1)^2\\ 
 &\le& 1.
\end{eqnarray*}
Therefore, $w(A_1)\le 1$. By the Assertion,  $|s_1| \le \hat s_1$.
Hence 
$A=\(1-\dfrac{|s_1|}{\hat s_1}\)A_0+\dfrac{|s_1|}{\hat s_1} A_1$ 
is a convex combination of $ A_0$ and $A_1$.

\medskip
Similarly, if  $W(B)$ touches the unit circle at $e^{i\phi_2}$
with $\phi_2\in[-\pi/2,\pi/2]$, then  $B$ is a 
convex combination of 
\begin{equation}\label{b0b1}
B_0=e^{i\phi_2}I \qquad  \hbox{ and } \qquad  B_1 = i(1-r^2)\sin\phi_2I+\nu_2\hat s_2 C
\end{equation}
with $\hat s_2 = \sqrt{\cos^2\phi_2+r^2 \sin^2 \phi_2}$ and $\nu_2\in\{1,-1\}$.  
Let $U = \begin{pmatrix}-r& \sqrt{1-r^2}   \cr    \sqrt{1-r^2}& r\end{pmatrix}$. Then 
$U^*CU=-C$. If $\nu_2 = -1$, we may replace $(A,B)$ by $(U^*AU,U^*BU)$ 
so that $(\nu_1,\nu_2)$ will change to $(-\nu_1, -\nu_2)$. So, we may further 
assume that $\nu_2=1$.

\medskip
By the above analysis, 
$AB$ is a convex combination of $A_0B_0, A_0B_1, A_1B_0$ and $A_1B_1$. Since
$w(e^{it}T)=w(T)$ for all $t\in\IR$ and $T\in M_n$, 
the first three matrices have numerical radius 1.
We will prove that 
\begin{equation}
\label{a1b1}
w(A_1 B_1) < 1.
\end{equation} 
It will then follow that 
$w(AB) \le 1$, where the equality holds only when $A=A_0$ or $B=B_0$.

 \medskip
For simplicity of notation,
let $w_1=\sin\phi_1$ and $w_2=\sin\phi_2$. Then 

\begin{equation}\label{hats}\hat s_i=\sqrt{1-(1-r^2)w_i^2}\ \mbox{ for }\ i=1,2.
\end{equation}

Recall  from (\ref{a0a1}) and (\ref{b0b1}) that  $A_1 =   i(1-r^2)w_1 I +\nu_1 \hat s_1 C$ and
$B_1 =  i(1-r^2)w_2 I +  \hat s_2 C$ because $\nu_2=1$.   Since $C^2 = (1-r^2)I_2$, we have 
$$A_1B_1 = (1-r^2)(u I_2 + iv C),$$
where 
$$
u = \nu_1 \hat s_1\hat s_2-w_1w_2(1-r^2) 
\quad\hbox{ and } \quad v = w_1 \hat s_2
+\nu_1w_2\hat s_1.
$$
If $r=1$, then $A_1B_1 =0$. Assume that $0<r<1$.
We need to show that 
$$
\frac{1}{1-r^2}w(A_1B_1) = w(uI+ivC) < \dfrac{1}{(1-r^2)}.
$$
Because $W(uI+ivC)$ is an elliptical disk with boundary
$\{u+iv(\cos\theta+ir\sin\theta): \theta \in [0, 2\pi]\}$,
it suffices to show that
%\begin{equation}\label{uv0}
$$f(\theta) = |u + iv(\cos\theta + i r\sin \theta)|^2
  < \frac{1}{(1-r^2)^2}
\quad \hbox{ for all } \ \theta \in [0, 2\pi].$$
%\end{equation}
Note that
\begin{eqnarray*}
f(\theta) &=& (u-rv\sin\theta)^2 + (v\cos\theta)^2 \\
&=&u^2-2ruv\sin\theta+r^2v^2\sin^2\theta+v^2(1-\sin^2\theta) \\
&=&\dfrac{u^2}{1-r^2}+v^2-\(\sqrt{1-r^2}v\sin\theta+\dfrac{ru}{\sqrt{1-r^2}}\)^2\\
&\le&\dfrac{u^2}{1-r^2}+v^2\\
&=&\dfrac{1}{(1-r^2)}\[u^2+(1-r^2)v^2\]\\
&=&\dfrac{1}{(1-r^2)}\[ (\nu_1 \hat s_1\hat s_2-w_1w_2(1-r^2) )^2+(1-r^2)(w_1 \hat s_2+\nu_1w_2\hat s_1)^2  \]\\
&=&\dfrac{1}{(1-r^2)}\[ \hat s_1^2\hat s_2^2+w_1^2w_2^2(1-r^2)^2 +(1-r^2)
(w_1^2 \hat s_2^2+ w_2^2\hat s_1^2)\]\quad 
\mbox{because }\nu_1=\pm 1 \\
&=&\dfrac{1}{(1-r^2)}\[ (\hat s_1^2+ (1-r^2) w_1^2)(  \hat s_2^2   +(1-r^2) w_2^2)\]\\
&=&\dfrac{1}{(1-r^2)}\hskip 3.25in \ \mbox{ by (\ref{hats})}   \\
&< &\dfrac{1}{(1-r^2)^2}\hskip  3.25in \ \mbox{because }0<r<1.
\end{eqnarray*}
Consequently, we have   $w(A_1B_1)<1$ as asserted in (\ref{a1b1}).
Moreover, by the comment after (\ref{a1b1}), if
$w(AB) = w(A)w(B)$, then $A = A_0$ or $B = B_0$.
Conversely, if $A = A_0$ or $B_0$, then we clearly have $W(AB) = w(A)w(B)$.
The proof of the theorem is complete. \qed

\bigskip\noindent
{\bf Acknowledgment}

We would like to thank Professor Pei Yuan Wu, Professor 
Hwa-Long Gau, and the referee for some helpful comments. 
Li is an affiliate member of the Institute
for Quantum Computing, University of Waterloo, and is an
honorary professor of the
Shanghai University. His research was supported by USA
NSF grant DMS 1331021, Simons Foundation Grant 351047,
and NNSF of China Grant 11571220.

\end{document}